\documentclass[12pt]{article}
\usepackage{graphicx}
\usepackage{amsmath}
\usepackage{amssymb}
\usepackage{amscd,color}

\newtheorem{Rem}{Remark}
\newtheorem{theorem}{Theorem}
\newtheorem{theorema}{Theorem}

\newtheorem{lem}[theorema]{Lemma}

\newenvironment{proof}[1][Proof]{\textbf{#1.} }{\qed \vspace{5pt}}

\newcommand\bib[1]{\bibitem[#1]{#1}}
\newcommand\qed{\phantom{\underline{y}}\hfill\hfill$\square$}
\newcommand{\comm}[1]{}

\renewcommand\a{\alpha}
\renewcommand\b{\beta}

\renewcommand\d{\delta}

\renewcommand\l{\lambda}

\newcommand\op[1]{\mathop{\rm #1}\nolimits}

\newcommand\R{{\mathbb R}}

\newcommand\we{\wedge}
\newcommand\x{\xi}

\def\R{\mathord{\mathbb R}}

\def\so{\mathfrak{so}}

\def\2{\frac{1}{2}}
\def\3{{\ss}}
\def\x{\times}
\def\.{\cdot}
\def\d{\partial}

\def\<{\langle}
\def\>{\rangle}

\def\Hom{\mathop{\rm Hom}\nolimits}

\def\ms{\medskip\noindent}

\def\beq{\begin{equation}}
\def\eeq{\end{equation}}
\def\bea{\begin{eqnarray}}
\def\eea{\end{eqnarray}}
\def\bsm{\left(\begin{smallmatrix}}
\def\esm{\end{smallmatrix}\right)}
\def\bpm{\begin{pmatrix}}
\def\epm{\end{pmatrix}}

\newcommand{\weg}[1]{}

\begin{document}

 \title{Compatibility of Gau\3 maps with metrics}
 
 \author{J.\ Eschenburg, B.\ S.\ Kruglikov, V.\ S.\ Matveev and 
 R.\ Tribuzy} \date{}
 \maketitle

 \begin{abstract}
We give necessary and sufficient conditions on a smooth local map of a Riemannian
manifold $M^m$ into the sphere $S^m$ to be the Gau\3 map of an isometric immersion
$u:M^m\to\R^n$, $n=m+1$. We briefly discuss the case of general $n$ as well\\[.1cm] 

 \end{abstract}

\section{Introduction}

Isometric hypersurface immersions of a Riemannian manifold $(M,g)$ with
dimension $m = n-1$ into Euclidean $n$-space are characterized by their first
and second fundamental forms, $g$ and $h$. By a classical theorem going back to
Bonnet, the immersion  exists and
is uniquely determined by $g$ and $h$ up to Euclidean
motions provided that the pair $(g,h)$ satisfies Gau\3 and Codazzi equations.

In the present paper we ask what happens if we replace $(g,h)$ by $(g,\nu)$
where $\nu : M \to S^m$ is the Gau\3 map. At the first glance the new problems
seems more rigid since $h$ is obtained from the differential $d\nu : TM \to
\nu^\perp$. However this observation is true only after identifying $TM$ with
the complement of the normal bundle, $\nu^\perp$. This identification is
precisely  the differential of the immersion which has to be constructed.
In fact uniqueness might fail as it happens with minimal surfaces:
All immersions in the associated family of a minimal surface have the
same $g$ and $\nu$, but they are not congruent.
A well-known example is the deformation of the catenoid to the helicoid.

The first problem is to recover the second fundamental form $h$ from the data.
In our approach, the third fundamental form $k = \<d\nu,d\nu\>$ will play a
major r\^ole since it is obtained directly from our data and the second
fundamental form is its square root (using $g$, all 2-forms are viewed also
as endomorphims).
However, the square root of a self adjoint positive
semi-definite matrix is not unique, and if repeated eigenvalues occur, there
are even infinitely many solutions as it happens in the minimal surface case
mentioned above.
Moreover, in high dimensions it might be very difficult to compute.

Fortunately, Theorems \ref{thm2}, \ref{thm3}
give  other ways to recover $h$  (using Gau\3 equations).
The defining equations
(\ref{h2}) and (\ref{h=}) in theorems 1 and 2 have been known already to Obata
\cite{O}. In the  final chapter we extend these ideas to the case of higher
codimension where $\nu$ takes values in the Grassmannian $G_{n-m}(\R^n)$.

\section{Main results} \label{s3}

Let $(M^m,g)$ be a Riemannian manifold and $\nu:M^m\to S^m\subset\R^{m+1}$  be { a smooth mapping}.
We are interested in the question when the given data
$(g,\nu)$ are the first fundamental form and the Gau\ss\ map
for an immersion $u:M^m\hookrightarrow\R^{m+1}$. Such data  $(g,\nu)$ will be
called {\it admissible}, and $u$ will be called a {\it solution} for $(g,\nu)$.
Our considerations are local, hence we may always assume that $M$
is a simply connected open subset of $\R^m$.

Let $R{ \,=\,(R^l_{ijk}) }$ be the Riemann curvature tensor of the metric $g$, $\op{Ric} =
\op{Tr} R{\, = \left( R^l_{ijl} \right)}$
the Ricci tensor,  and $s = \op{Tr}\op{Ric}= g^{ij}\op{Ric}_{ij}$ the scalar curvature.
Let $A = d\nu \in \Hom(TM,\nu^\perp)$ and put
\beq \label{k}
	k(v,w)=\langle Av,Aw\rangle = \<A^*Av,w\>
\eeq
for all $v,w\in TM$;  this is a symmetric positive semi-definite bilinear
form, which will be referred to as the \emph{third fundamental form}.
We can raise the indices with the help of $g$ and consider both $\op{Ric}$ and  $k$
as fields of operators on the tangent bundle, denoting the result by the same
letter.

 \begin{theorem}\label{thm1}
Let $(M^m,g)$ and $\nu : M^m \to S^m$ be given and assume that $A = d\nu : TM^m \to
\nu^\perp$ is
everywhere invertible. Let $k$  be defined by
{\rm (\ref{k})}.
Then the data $(g,\nu)$ are admissible
if and only if there is $h \in S^2T^*M$ with
\beq \label{h2}
	h^2 = k.
\eeq
such that the vector bundle homomorphism
  \beq \label{U}
U = -(A^*)^{-1}h : TM \to \nu^\perp
\eeq
is isometric 
and parallel with respect to the Levi-Civita connection on $TM$
and the projection connection on $\nu^\perp$.
In fact, the corresponding immersion $u : M^m \to \R^{m+1}$ is determined by
$\,du = U$, and $h$ is the second fundamental form of $u$, i.e. $h_{ij} = \<u_{ij},\nu\>$.
\end{theorem}

Since $k$ is positive semi-definite, it has a symmetric square root $h$.
However, as explained in the introduction,
(\ref{h2}) is often difficult to solve. Indeed,   a  general solution requires finding  the roots of  a polynomial  of the $m^{\textrm{th}}$ degree. For big $m$, this  is  impossible to do  explicitly.
If $k$ has multiple eigenvalues, the following additional difficulty appears:   at every point
there are infinitely many  solutions of (\ref{h2}), so even if we found one solution of (\ref{h2}) such that \eqref{U} is not parallel, there
might exist another  solution such that \eqref{U} {\emph is}  parallel
Hence in many cases  Theorem \ref{thm1}  is useless unless we find a better
method to compute $h$ from the data. This is achieved by the following statements:

 \begin{theorem} \label{thm2} 
If the data $(g,\nu)$
are admissible, then for any solution $u$, the second fundamental form $h$ and the (unnormalized) mean curvature $H = \op{Tr} h = h_{ij}g^{ij}$ solve the following system
  \bea\label{h=}
h\,H = \op{Ric}+k,\quad  H^2 = s+\op{Tr}k.
  \eea
 \end{theorem}

\begin{Rem} 
Clearly, if $s+ \op{Tr}k> 0$, the equations can be solved: 

 \bea\label{h=1}
  H = \pm\sqrt{ s+\op{Tr}k}, \quad h=\pm \tfrac{1}{\sqrt{s+\op{Tr}k} } \left(\op{Ric}+k\right).
  \eea
  Moreover,  as we explain in Remark \ref{rem3}, the sigh of $H$ and of $h$ is  not essential for our goals.   
\end{Rem} 

\begin{theorem}
\label{thm3} If the data $(g,\nu)$
are admissible with $d\nu$ non-degenerate and $m=\dim M\ge 3$, then the second fundamental form $h$ of any solution $u$  solves the homogeneous linear system
\beq \label{lin}
	h\,k^{-1}R(\Omega) = 2\,\Omega\,h,\quad \forall\, \Omega \in \so_g(TM),
\eeq
where $R$ is considered as curvature operator acting on $\Lambda^2TM = \so_g(TM)$. 
Moreover the solution $h$  of \eqref{lin} is unique up to a scalar factor.
\end{theorem}

\begin{Rem} \label{rem2}
Note that the missing scalar factor in Theorem \ref{thm3} can be easily found using 
condition \eqref{h2}: if $\tilde h$ is a nonzero  solution of \eqref{lin} then
 \beq \label{scalar}
h= \pm \sqrt{\frac{ \op{Tr} \left( \tilde h^2 \right)}{\op{Tr} k}} \cdot \tilde h
 \eeq
 \end{Rem} 
 
 \begin{Rem} \label{rem3} 
The sign $\pm$ in the formulas (\ref{h=1},\ref{scalar}) does not affect the existence of a solution $u$: 
If $u$ is a solution for $(g,\nu)$ with second fundamental
form $h$ (resp. mean curvatrure $H$), then $-u$ is also  a solution with second fundamental form $-h$ (resp. mean curvature $-H$).
\end{Rem}

\noindent The above theorems give us an algorithm to check admissibility of $(g,\nu)$:

\begin{enumerate}
\item Check if   $ s+ \op{Tr}k \geq 0$.\footnote 
	{In fact, a bit more is necessary: $s+ \op{Tr}k$ needs to allow a smooth
	``square root'': a function $H$ with $H^2 = s+ \op{Tr}k$. }
\item  Find $h$:
\begin{enumerate}
\item If  $s+ \op{Tr}k> 0$, define $h$ by (\ref{h=1}).\footnote
	{The sign of $H$ is arbitrary, see Remark \ref{rem3}.}
 \item If $s+ \op{Tr}k=0$ 
and $m \geq 3$, then  for every $x\in M$  solve the  linear system (\ref{lin}) in $T_xM:= \nu^\perp$. Check whether there exists a  (nondegenerate)  solution $\tilde h$.   Consider the solution $h$ given by \eqref{scalar}.
 \item If $s+ \op{Tr}k\equiv 0$ and $m=2$, verify the Gau\ss\ condition (Remark \ref{rem6} below).
  \end{enumerate}
\item Check if $h^2 = k$ (this together with 2 implies that $h$ is symmetric). 
\item Finally check if $U = - (A^*)^{-1}h$
is parallel, i.e.
$$
	\d_iu_j - \<\d_iu_j,\nu\>\nu = \Gamma_{ij}^k u_k
$$
where $u_j = Ue_j$ and $\Gamma_{ij}^k$ are the Christoffel symbols, the
components of the Levi-Civita connection: $\nabla_ie_j = \Gamma_{ij}^ke_k$.
\end{enumerate}
The data $(g,\nu)$ are admissible if and only if  all checks are successful.

This answers a question raised in \cite{E} which was discussed by the authors
during and after 
the 10$^\text{th}$ conference on Differential Geometry and its
Applications.

\begin{Rem} 
 Both Gau\ss\ and Codazzi equations are hidden in the assumption
that $U \in \Hom(TM,\nu^\perp)$ is isometric and 
parallel. In fact this property
is equivalent to Codazzi equations while Gau\3 equations follow from it,
see Appendix. The claim that the Gau\ss\ condition is a differential
consequence of the Codazzi condition in the non-degenerate case is
non-trivial. It shall be compared with the known
fact that under some conditions the Codazzi equations are consequences of the
Gau\ss\ equations \cite{Al}. \end{Rem}

\begin{Rem} The uniqueness of recovering the isometric immersion
$u:M^m\to\R^{m+1}$ with fixed third quadratic form $k$ was considered in
\cite{DG}. This is similar  to recovering of immersions with fixed
Gau\ss\ map in the case of hypersurfaces, but not for higher codimension,
see the last section. \end{Rem} 

\begin{Rem} \label{rem6}  The only case not covered by our theorems is $m = 2$ and $H = 0$,
the case of minimal surfaces which is given by the well known Weierstra\3
representation \cite{L}; the only restriction for the metric $g$ comes from the Gau\3
equation
 $$
K+\sqrt{\op{det}(k)}=0
 $$
and the Gau\3 map $\nu : M^2 \to S^2$ needs to be conformal. Any such pair
$(g,\nu)$ is admissible, and to each admissible pair there exists precisely
a one-parameter family of geometrically distinct isometric minimal immersions,
the associated family. \end{Rem}

\section{Historical motivation}

Two classical problems concern the embeddings
 \begin{equation}\label{E0}
u:M^m{ \hookrightarrow}\R^{n}.
 \end{equation}

The first is about isometric embedding, i.e. when a metric $g$ on $M$ can be obtained as $u^*ds^2_\text{Eucl}$ for some $u$. In the PDE language this is equivalent to solvability
of the system
 \begin{equation}\label{E1}
\langle u_{i},u_{j}\rangle=g_{ij}(x),\qquad 1\le i,j\le m,\quad
u=(u^1,\dots,u^n):M^m\to \R^n,
 \end{equation} where $u_i:= \tfrac{\partial }{\partial x_i}u$.

The Janet-Cartan theorem \cite{J,C} guarantees this locally in the analytic category
for $n=\frac{m(m+1)}2$, i.e. when the system (\ref{E1}) is determined. This was
improved by Nash \cite{N}, Gromov-Rokhlin \cite{GR} and others \cite{Gre,BBG},
who relaxed analyticity to smoothness (for the price 
of increasing $n$ or imposing some non-degeneracy assumptions) and so proving that
embedding is always possible. When $n<\frac{m(m+1)}2$ the system is
overdetermined. Thus while rigidity (uniqueness of solutions up to Euclidean
motion) is clear in many cases, no general criterion (existence) for local
embedding is known (see \cite{Gro} for details).

The other important problem related to imbeddings (\ref{E0}) is to recover it from the Gau{\ss} map $\nu:M^m\to G_{n-m}(\R^n)$, $x\mapsto T_xM^\perp\subset\R^n$ (also known as Grassmann map). This problem is unsolvable for hypersurfaces ($n=m+1$) unless the Gau\ss\ map is degenerate. In general the problem can be rephrased as solvability of the following PDE system
 \begin{equation}\label{E2}
\langle u_{i},\nu^\alpha(x)\rangle=0,\quad m+1\le \alpha\le n
 \end{equation}
where $\nu^\alpha$ is an orthonormal basis of sections of $TM^\perp$
(no index for hypersurfaces).

For $2m=n=4$, the system (\ref{E2}) is determined while for
the other $m<n-1$ overdetermined. By the results of Muto, Aminov,
Borisenko \cite{M,Am,B$_1$} the embedding is locally recoverable upon
certain non-degeneracy assumptions, up to a parallel translation and homothety.
However not any $m$-dimensional submanifold of $G_{n-m}(\R^n)$ is realizable as the image of a Gau\ss\ map (except for the case 2 in 4,
when no obstruction equalities exist). The conditions of realizability are not
known so far (partial results can be found in \cite{B$_2$}).

In this note we unite the systems (\ref{E1})+(\ref{E2}) and ask when the data $(g,\nu)$
are realizable and what is the freedom. In many cases we get indeed rigidity, i.e. an embedding is recoverable up to a parallel translation (this can be obtained as a
combination of the problems with the data $g$ and the data $\nu$ above, but
our conditions are wider; another approach to rigidity within the same
problem was taken in \cite{AE}). However in addition to this we write the
full set of constraints, thus solving the problem completely.

Notice that as $(n-m)$ grows, the amount of compatibility constraints coming from $g$ decreases while that for $\nu$ increases, and there are always constraints for $(g,\nu)$.
For the hypersurface 
case $n=m+1$, mainly treated here, the Gau\ss\ map alone
bears no information (unless it is degenerate),
making the problem for the pair $(g,\nu)$ more interesting.

\section{Proof of the main results}\label{S3}
We will work locally in $M^m$. Given $(M,g)$ and $\nu: M^m \to S^m$,
we want to understand whether there exists  a smooth (local) map  $u : M^m \to \R^{ m +1}$
whose partial derivatives with respect to a coordinate chart
satisfy
\bea	\label{uiuj}
	\<u_i,u_j\> &=& g_{ij},	\\
	\label{uinu}
	\<u_i,\nu\> &=& 0.
\eea
This is a  system of algebraic equations for the partial derivatives of $u_i$.
Once we obtained
an (algebraic) solution $U = (u_i)_{i=1,\dots,m}$ of
(\ref{uiuj}),(\ref{uinu}), {  there exists  a smooth mapping  $u= (u^1,...,u^n):
M \to \R^n$  }
if and only if the
integrability conditions
\beq	\label{uij}
	u_{ij} = u_{ji}
\eeq
are fulfilled, where the second index
means partial derivative: { $u_{ij}:= \tfrac{\partial}{\partial x_j}u_i$}.
Equation (\ref{uij}) splits into a tangent and a normal
part. The tangent part ($\nu^\perp$-part) can be
interpreted as follows. Equations (\ref{uiuj}), (\ref{uinu})
mean that $U$ is a
bundle isometry between $TM$ and $\nu^\perp$. Now the tangent
part of (\ref{uij})
says that the canonical connection (via projection) on
$\nu^\perp \subset M \x \R^n$ is torsion free when
$\nu^\perp$ and $TM$ are identified using $U$. Since the connection
is also metric preserving, it is Levi-Civita:
\beq \label{uijT}
	(u_{ij})^T = \Gamma_{ij}^k u_k
\eeq
where $(\ )^T$ denotes the tangent component ($\nu^\perp$-component).
In other words, the tangent part of (\ref{uij}) under the assumptions
(\ref{uiuj}),(\ref{uinu}) says
precisely that $U : TM \to \nu^\perp$ is parallel (affine, connection
preserving).

\ms The normal part ($\nu$-part) of (\ref{uij}), in view of (\ref{uinu}), is
equivalent to
 $$
h_{ij} = h_{ji}
 $$
where
\beq \label{h=AU}
	h_{ij} = \<u_{ij},\nu\> = -\<u_i,\nu_j\>,\ \ \ h = -A^*U.
\eeq
Once we have got $h$, we obtain $U = -(A^*)^{-1}h : TM \to \nu^\perp$ from
(\ref{h=AU}) and check orthogonality (\ref{uiuj}) and parallelity (\ref{uijT}).

Next we show that $h^2=k$ is necessary. If an immersion $u : M \to \R^n$ with
Gau\3 map $\nu$ is given, then $h = -A^*U$ where $A = d\nu$ and $U = du$ because
$h_{ij} = \<u_{ij},\nu\> = -\<u_i,\nu_j\>$.
Since $h$ is self adjoint and $U$ orthogonal, we have
$$
	h^2 = hh^* = A^*UU^*A = A^*A = k.
$$

Now let us show that our assumptions are sufficient. Assuming
$h$ symmetric with $h^2=k$ 
and choosing $U = -(A^*)^{-1}h : TM \to \nu^\perp$, we obtain
 $$
UU^* = (A^*)^{-1}h^2A^{-1} = (A^*)^{-1}kA^{-1} = (A^*)^{-1}A^*AA^{-1} = \textrm{\rm I},
 $$
thus $U$ is an isometry. Moreover $h = -A^*U$ is symmetric, i.e.
 $$
\<u_i,\nu_j\> = \<u_j,\nu_i\>.
 $$
Since $U$ takes values in $\nu^\perp$, we have $\<u_i,\nu\> = 0$ and hence
 $$
\<u_{ij},\nu\> = \<u_{ji},\nu\>.
 $$
This is the normal part of (\ref{uij}). The tangent part is obtained from
the parallelity assumption (\ref{uijT}), since the Christoffel
symbols $\Gamma_{ij}^k$ are symmetric in $(ij)$. Thus the
integrability condition (\ref{uij}) is proved and hence we obtain a map $u : M \to \R^n$
with $du = U$. This finishes the proof of Theorem \ref{thm1}.

\medskip

Theorem 2 is an obvious consequence of the Gauss equations
\beq \label{gauss}
	R_{ijkl} = h_{il}h_{jk} - h_{ik}h_{jl}.
\eeq
In fact, taking the trace over $jk$, i.e. multiplying by $g^{jk}$ and summing
we obtain
\beq \label{Ric}
	\op{Ric} = h\.H - k,
\eeq
Taking again the trace on both sides,
 \beq \label{kH}
	s = H^2 - \op{Tr}k
 \eeq
This shows (\ref{h=}). Theorem \ref{thm2} is proved.

In order to prove Theorem \ref{thm3},
 we transform equations (\ref{gauss})
into its curvature operator form
$$
R(\Omega) = 2\. h\Omega h,    \  \  \  \textrm{{ i.e., { ${R^{ij}}_{kl} \Omega^{kl} =  2 \. h^i_l \Omega^{kl} h^j_k$}}}
 $$
(which must be fulfilled for every   $\Omega \in \Lambda^2TM = \so_g(TM)$).
Multiplying by $h^{-2} =
k^{-1}$ from the left, we get
\beq \label{Adh}
	k^{-1}R(\Omega) =2\cdot  h^{-1}\Omega h
\eeq
which is equivalent to
\beq \label{lin2}
	hk^{-1}R(\Omega) = 2\cdot \Omega h.
\eeq
This is the linear equation (\ref{lin})  for $h$ which we wanted to prove. It remains
to show  uniqueness of the solution provided $m\ge 3$. This will be done in the following

\begin{lem} Assume $m\ge 3$.
If both $h,\tilde h$ solve {\rm (\ref{lin2})} for all $\Omega\in\so(\mathbb{R}^m)$,
then $h$ and $\tilde h$ are proportional. 
\end{lem}

\noindent
\proof Both $h$ and $\tilde h$ satisfy (\ref{Adh}) and hence
$$
	h^{-1}\Omega h = \tilde h^{-1}\Omega \tilde h
$$
for all $\Omega \in \so(\mathbb{R}^m )$. Thus $g := \tilde h h^{-1}$ commutes with all
$\Omega \in \so(\mathbb{R}^m )$. If $m\geq 3$, the centralizer of $\so(\mathbb{R}^m)$ contains only the
scalar matrices, hence $g = \lambda I$ and $\tilde h = \lambda h$ for some
$\lambda \in\R$.
\endproof

\section{The general case}

When we study immersions $M^m\to\R^n$ with general $n$, it is not
easy to get a closed formula for $h$ as in the main theorems.
But as we will see, in the generic case $h$ can still 
be effectively computed from the data.

The given data are again a Riemannian manifold $(M^m,g)$ and a smooth map $\nu
: M \to G_{n-m}(\R^n)$ into the Grassmannian of $(n-m)$-planes in $\R^n$, and
we ask if $\nu$ is the Gau\ss\ map of some isometric immersion $u : M \to \R^n$.
Choose an orthonormal basis $(\nu^\alpha)_{\alpha = 1,\dots,n-m}$ of $\nu$.
Let
 $$
A^\alpha = (d\nu^\alpha)^T
 $$
denote the corresponding Weingarten operators, where $(\ )^T$ again denotes the
tangent component ($\nu^\perp$-component), and let
\beq \label{kalphabeta}
	k = \sum k^{\alpha\alpha},\ \ \ k^{\alpha\beta} = (A^\alpha)^* A^\beta
\eeq
the corresponding third fundamental on $M$ induced by $\nu$
from the standard symmetric metric on $G_{n-m}(\R^n)$.
The second fundamental form $h$ which we
search for, is $\nu$-valued and has also several components $h^\alpha =
\<h,\nu^\alpha\>$: Given an isometric immersion $u : M \to \R^n$ and $U =du=
(u_1,\dots,u_m)$, we have
 $$
h^\alpha_{ij} = \<u_{ij},\nu^\alpha\> = - \<u_i,\nu^\alpha_j\>,\qquad h^\alpha = -(A^\alpha)^*U.
 $$
Consequently
 $$
h^\alpha h^\beta = h^\alpha(h^\beta)^* = (A^\alpha)^* UU^*A^\beta	= k^{\alpha\beta}.
 $$

There exist deformations of isometric immersions with fixed $ k$, see \cite{V},
which are different from ours for codimensions exceeding 1. However in this case
$ k$ bears significantly less information than the Gau\ss\ map $\nu$. With the
latter we can restore the immersion up to a translation in a generic case.

 \begin{theorem}
Let $(M^m,g)$ and $\nu : M \to G_{n-m}(\R^n)$ be given and assume that
at least one of the $A^\alpha = (d\nu^\alpha)^T : TM \to \nu^\perp$ is
everywhere invertible. Let $k^{\alpha\beta}$ be as in
{\rm (\ref{kalphabeta})}.
Then the data $(g,\nu)$ are admissible
if and only if there exist $h^\alpha \in S^2T^*M$ with
\beq \label{halpha2}
	h^\alpha h^\beta = k^{\alpha\beta}.
\eeq
and a vector bundle homomorphism $U : TM \to \nu^\perp$ with
\beq \label{AalphaU}
	h^\alpha = -(A^\alpha)^*U
\eeq
for all $\alpha$, such that $U$ is parallel with respect to the
Levi-Civita connection on $TM$
and the projection connection on $\nu^\perp$.
In fact, the corresponding immersion $u : M^m \to \R^n$ is determined by
$du = U$.
 \end{theorem}

The proof is almost the same as before and will be omitted.
But as before we need an effective method to compute
$h^\alpha$ from the given data. This is given by the next
theorem:

\begin{theorem} Assume that the data $(g,\nu)$
are admissible with $|H|=\sqrt{s+\op{Tr}k}$ $\neq 0$ and $\op{Ric}+k$
invertible. Then
\beq \label{hbeta=}
	h^\beta =\sum_\alpha H^\alpha(\op{Ric}+k)^{-1}k^{\alpha\beta}.
\eeq
where $H^\alpha = \op{Tr}h^\alpha$ are the components of the mean curvature
vector $H = \op{Tr}h$ which is a fixed vector with length $\sqrt{s+\op{Tr}k}$ 
for the matrix $\rho = (\rho_{\alpha\beta})$ on $\nu$ defined by 
\beq \label{rho}
	\rho_{\alpha\beta} = \op{Tr}\left((\op{Ric}+k)^{-1}\right)k^{\alpha\beta}.
\eeq
\end{theorem}

\noindent
\proof
Suppose that an isometric immersion $u : M \to \R^n$ with Gau\3 map $\nu$ is
given. The Gau\3 equations are
 $$
R_{ijkl} = \sum_\alpha h^\alpha_{il}h^\alpha_{jk} - h^\alpha_{ik}h^\alpha_{jl}
 $$
Taking the trace over $jk$ yields
$$
	\op{Ric} = \sum_\alpha \left(h^\alpha H^\alpha - (h^\alpha)^2\right)
	= \sum_\alpha h^\alpha H^\alpha - k
$$
and hence
\beq \label{Hh}
	\sum_\alpha H^\alpha h^\alpha = \op{Ric} + k.
\eeq
Tracing again we obtain the length of the mean curvature vector,
\beq \label{|H|}
	|H| = \sqrt{s+\op{Tr}k}.
\eeq
From (\ref{Hh}) we can compute the $h^\alpha$ since the products $h^\alpha
h^\beta = k^{\alpha\beta}$ are known:
$$
	(\op{Ric} + k)h^\beta = \sum_\alpha H^\alpha k^{\alpha\beta}
$$
and (\ref{hbeta=}) follows.
In order to compute $H^\alpha$ we take the trace of (\ref{hbeta=}):
 $$
	H^\beta =
\sum_\alpha H^\alpha \op{Tr}\left((\op{Ric} + k)^{-1} k^{\alpha\beta}\right)
= \sum_\alpha H^\alpha \rho_{\alpha\beta}
 $$
with $\rho_{\alpha\beta}$ as in (\ref{rho}). Thus $H$ is a fixed vector of
the matrix $\rho = (\rho_{\alpha\beta})$. In the generic case this fixed space
is only one dimensional (it is at least one-dimensional since it contains $H \neq
0$). Using (\ref{|H|}) we see that $H$ is uniquely determined up to
sign.\footnote
	{The sign of $H$ cannot be fixed. Indeed, as in Remark \ref{rem3},  if $u : M \to \R^n$ is
	an immerion with mean curvature vector $H$, then $-u$ has the same
	Gau\3 map and mean curvature vector $-H$.}
\endproof

\noindent Once again we have got an algorithm by which we may check 
if the data $(g,\nu)$ are
admissible, belonging to some isometric immersion $u$. From the data we form the
matrix $\rho$ and check if it has a fixed vector. In the generic case, the fixed
space is at most one-dimensional. We choose a fixed vector $H$ using
(\ref{|H|}). Then we define the quadratic form $h^\beta$ by (\ref{hbeta=}) and
check if it satisfies (\ref{halpha2}). Setting $U =-(A^{\alpha*})^{-1}h^\alpha$
for one $\alpha$ we check if (\ref{AalphaU}) holds true for the other indices $\alpha$ and
if $U$ is parallel. The data $(g,\nu)$ are admissible if and only if all the tests
are successful.

In non-degenerate cases, similar to the hypersurface case, the Gau\ss\ equation follows from Codazzi and Ricci equations.
Non-uniqueness of solution for this system means isometric deformation with
fixed Gau\ss\ image, and can be considered similarly to what was done for
$n=m+1$. We expect that the only examples are locally products of the
(cylinders over) Weierstra\ss\ examples, so that an essentially 2-dimensional
phenomenon generates all such examples.

\section*{Appendix}

We want to show that under the assumptions (\ref{uiuj}), (\ref{uinu}) and
(\ref{uij})$^\perp$, Codazzi equations are equivalent to
parallelity of $U$, (\ref{uij})$^T$, 
and they imply Gau\ss\ equations.

 \begin{theorem}
Let $(M,g)$ and $\nu : M \to S^m$ be given and $d\nu$ nondegenerate. Let
$U= (u_1,\dots,u_m) : TM \to \nu^\perp$
be a vector bundle isometry such that $b_{ij}: = \<u_i,\nu_j\> = -h_{ij}$ is
symmetric (normal integrability condition). Then
\bea \label{codazzi}
\nabla U = 0 &\iff&
(\nabla_ib)_{jk}=(\nabla_{j}b)_{ik},	\\
\nabla U = 0 &\Rightarrow& R_{ijkl} = b_{il}b_{jk} - b_{ik}b_{jl}
\label{gaussgl}
\eea
\end{theorem}

\noindent
 \begin{proof} ``$\Rightarrow$'' of (\ref{codazzi}):
From $b = U^*d\nu = \<U,d\nu\>$ we obtain
\beq \label{nablab}
	\nabla b =  \<\nabla U,d\nu\> + \<U,\nabla d\nu\>.
\eeq
Since $\nabla U = 0$, we are left with $\nabla b = \<U,\nabla d\nu\>$ or more precisely,
$$
	(\nabla_ib)_{jk} = \<u_j,(\nabla_i d\nu)_k\>
$$
Since the right hand side (the hessian of the map $\nu$) is symmetric in $ik$, we have
proved our claim.

\ms
``$\Leftarrow$'' of (\ref{codazzi}):
We still have (\ref{nablab}), more precisely
$$
	(\nabla_i b)_{jk} = \<(\nabla_i U)_j,\nu_k\> + \<u_j,(\nabla_i d\nu)_k\>.
$$
From the symmetry of $\nabla b$ and $\nabla d\nu$ in $ik$ we see
\beq \label{nablaUnu}
	\<(\nabla_i U)_j,\nu_k\> = \<(\nabla_k U)_j,\nu_i\>.
\eeq
On the other hand, by covariant differentiation of the isometry property $\<U,U\> = U^*U = g$
we obtain $\<\nabla U,U\> + \<U,\nabla U\> = 0$, more precisely
\beq \label{nablaUU}
	\<(\nabla_i U)_j,u_k\> + \<u_j,(\nabla_i U)_k\> = 0.
\eeq
Since $U^*d\nu = b$ is self adjoint with respect to $g$, we can choose local coordinates in such
a way that both tensors $g$ and $b$ are diagonal at the considered point and hence
$u_k=\lambda_k\nu_k$ for each $k$ (where we have used the nondegeneracy of $b$).
Substituting this into (\ref{nablaUU}) and putting
$\theta_{ijk} = \<(\nabla_i U)_j,\nu_k\>$,  we get
$$
	\lambda_k\theta_{ijk} + \lambda_j\theta_{jik} = 0
$$
(no summation). Cycling $(ijk)$ and using (\ref{nablaUnu}), $\theta_{ijk} = \theta_{kji}$,
we get 3 equations in the 3 unknowns $\theta_{ijk},\theta_{jki},\theta_{kij}$.
The determinant of this linear system
equals $2\l_i\l_j\l_k\ne0$, and therefore $\<(\nabla_i U)_j,\nu_k\>=\theta_{ijk}=0$.
Since the vectors $\nu_k$ form a basis of $\nu^\perp$, we obtain $\nabla U = 0$.

\ms Proof of (\ref{gaussgl}): Since $U : TM \to \nu^\perp$ is isometric and parallel, it
carries the Riemannian curvature tensor on $TM$ into the curvature tensor of
the projection connection $\nabla$ on $\nu^\perp$ which is computed as usual:
$$\begin{matrix}
	\nabla_j u_k &=& (u_{kj})^T &=& u_{kj} + b_{kj}\,\nu	\cr
	\nabla_i\nabla_j u_k &=& (\d_i(\nabla_j u_k))^T
	&=& (u_{kji})^T + b_{kj}\,\nu_i,	\cr
	\<\nabla_i\nabla_j u_k,u_l\> &=& \<\d_i(\nabla_j u_k),u_l\>
	&=& \<u_{kji},u_l\> + b_{kj}b_{il},	\cr
	R_{ijkl} &=&	\<(\nabla_i\nabla_j-\nabla_j\nabla_i)u_k,u_l\> &=&
	b_{kj} b_{il} -
	b_{ki} b_{jl},	
\end{matrix}$$
using the symmetry of $u_{kji} = \d_i\d_j u_k$ in $ij$.
The last equality is the Gau\3 equation (\ref{gaussgl}) which finishes the proof.
\end{proof}


\end{document}